\documentclass{amsart}

\usepackage{fancyhdr}
\usepackage{lastpage}
\usepackage{stmaryrd,yhmath}

\newcommand{\bdis}{\begin{displaymath}}
\newcommand{\edis}{\end{displaymath}}
\newcommand{\be}{\begin{equation}}
\newcommand{\ee}{\end{equation}}
\newcommand{\mbb}{\mathbb}
\newcommand{\mcal}{\mathcal}

\newcommand{\vp}{\varphi}

\newcommand{\vth}{\vartheta}

\newcommand{\zf}{\zeta\left(\frac{1}{2}+it\right)}

\theoremstyle{definition}

\theoremstyle{remark}
\newtheorem{remark}[]{Remark}

\newtheorem*{mydef1}{{\bf Theorem}}

\newtheorem*{mydef41}{{\bf Corollary 1}}

\newtheorem*{mydef42}{{\bf Corollary 2}}

\newtheorem*{mydef43}{{\bf Corollary 3}}

\newtheorem*{mydef7}{{\bf Question}}

\newtheorem*{mydef91}{{\bf Formula1}}

\newtheorem*{mydef92}{{\bf Formula2}}

\newtheorem*{mydef93}{{\bf Formula3}}

\newtheorem*{mydef94}{{\bf Formula4}}

\numberwithin{equation}{section}

\begin{document}

\title{Jacob's ladders, factorization and metamorphoses as an appendix to the Riemann functional equation for $\zeta(s)$ on
the critical line}

\author{Jan Moser}

\address{Department of Mathematical Analysis and Numerical Mathematics, Comenius University, Mlynska Dolina M105, 842 48 Bratislava, SLOVAKIA}

\email{jan.mozer@fmph.uniba.sk}


\begin{abstract}
In this paper we obtain a new set of metamorphoses of the oscillating Q-system by using the Euler's integral. We split the final state
of mentioned metamorphoses into three distinct parts: the signal, the noise and finally appropriate error term. We have also proved
that the set of distinct metamorphoses of that class is infinite one.
\end{abstract}
\maketitle

\section{Introduction and the first result}

\subsection{}

Let us remind the Riemann's functional equation (1859)
\be \label{1.1}
\zeta(1-s)=2(2\pi)^{-s}\cos\frac{\pi s}{2}\Gamma(s)\zeta(s),\ s\in\mbb{C}\setminus \{ 1\}.
\ee

\begin{remark}
It is known that L. Euler discovered a formula equivalent to (\ref{1.1}) for real values of the variable $s$ in
1749, (comp. \cite{2}, pp. 23-26), however the proof in the Euler's work is missing.
\end{remark}

Next, if we put
\be \label{1.2}
\begin{split}
& \chi(s)=\pi^{s-1/2}\frac{\Gamma\left(\frac 12-\frac 12s\right)}{\Gamma\left(\frac 12 s\right)} \ \Rightarrow \ \chi(s)\chi(1-s)=1,
\end{split}
\ee
(comp. \cite{9}, pp. 13, 16) then we obtain (see (\ref{1.1}), (\ref{1.2})) that
\be \label{1.3}
\frac{\zeta(1-s)}{\zeta(s)}=\chi(s).
\ee

\begin{remark}
Now, in connection with (\ref{1.3}) we use the following (E. Landau, \cite{3}, p. 30): the quotient
\bdis
\frac{\zeta(1-s)}{\zeta(s)}
\edis
is expressed by the \emph{known} function \dots (of course, \emph{known} function is such one that is different from $\zeta(s)$).
\end{remark}

\subsection{}

However, in the case
\be \label{1.4}
s=\frac 12+it,\ 1-s=\frac 12-it=\left(\frac 12+it\right)^\star
\ee
we have (comp. (\ref{1.2}), (\ref{1.4})) that
\bdis
\left| \chi\left(\frac 12+it\right)\right|=1,
\edis
i.e. (see (\ref{1.3}))
\be \label{1.5}
\left|\frac
{\overline{\zf}}{\zf}
\right|=1,\ \forall t\not=\gamma:\ \zeta\left(\frac 12+i\gamma\right)=0,
\ee
(after continuation for $t=\gamma$ this is valid for all $t\in\mbb{R}$).

\begin{remark}
Consequently, on the critical line
\bdis
s=\frac 12+it
\edis
the Riemann's functional equation (\ref{1.3}) gives only the trivial result (\ref{1.5}) (of course, the main aim of
the Riemann's functional equation is the analytic continuation of $\zeta(s)$ to $\mbb{C}\setminus \{ 1\}$).
\end{remark}

\subsection{}

It is clear that in this situation our Remark 2 leads us to the following

\begin{mydef7}
about another method of sampling of the points
\be \label{1.6}
t_2>t_1>0
\ee
(say) from a corresponding set that the quotient (comp. (\ref{1.3}), (\ref{1.5}))
\be \label{1.7}
\left|
\frac
{\zeta\left(\frac 12+it_2\right)}
{\zeta\left(\frac 12+it_2\right)}
\right|,\ t_1,t_2\not=\gamma
\ee
is expressed by a known function.
\end{mydef7}

There is method giving an answer to the Question -- namely one answer from a set of many possibilities -- mentioned method is our
method of transformation by using the reversely iterated integrals (comp. \cite{7}, (4.1) -- (4.19)).

\subsection{}

In this paper we obtain, for example, the following

\begin{mydef91}
For every sufficiently big $L\in\mbb{N}$ and for every $U\in (0,\pi)$ there are functions
\be \label{1.8}
\begin{split}
& \alpha_0^4=\alpha_0^4(L,U;a,b), \\
& \alpha_1^4=\alpha_1^4(L,U;a,b), \\
& \beta_1^4=\beta_1^4(L,U) , \\
& \alpha_0^4,\alpha_1^4,\beta_1^4\not=\gamma:\ \zeta\left(\frac 12+i\gamma\right)=0
\end{split}
\ee
such that
\be \label{1.9}
\begin{split}
& \left|
\frac
{\zeta\left(\frac 12+i\alpha_1^4\right)}
{\zeta\left(\frac 12+i\beta_1^4\right)}
\right|\sim \\
&
\sim
\frac{\arctan\left(\sqrt{\frac{a-b}{a+b}}\tan\frac U2\right)}
{\sqrt{\frac{a-b}{a+b}}\frac U2}
\frac{a+b\cos(\alpha_0^4)}{a+b},\ L\to\infty,
\end{split}
\ee
where
\be \label{1.10}
\begin{split}
& \alpha_0^4\in (2\pi L,2\pi L+U), \\
& \alpha_1^4,\beta_1^4\in (\overset{1}{\wideparen{2\pi L}},\overset{1}{\wideparen{2\pi L+U}}), \\
& (2\pi L,2\pi L+U)\prec (\overset{1}{\wideparen{2\pi L}},\overset{1}{\wideparen{2\pi L+U}}) ,
\end{split}
\ee
and the functions (\ref{1.8}) have the following property
\bdis
\alpha_1^4-\alpha_0^4\sim (1-c)\boldsymbol{ \pi}(2\pi L),\ L\to\infty,
\edis
where $c$ is the Euler's constant and $\boldsymbol{\pi}(t)$ is the prime-counting function.
\end{mydef91}

\begin{remark}
The formula (\ref{1.9}) gives one answer to our Question (comp. (1.7), (1.10)) in the direction outlined in Remark 2.
\end{remark}

\begin{remark}
Let us notice explicitly that in the case (\ref{1.7}), i.e. on the critical line, we may suppose that the \emph{known} function is \emph{every}
function, for example
\bdis
f[|\zeta|],\ \arg\zeta,\ \dots
\edis
\end{remark}

In this direction, we have obtained the following (see also corresponding formulae in the papers \cite{6} -- \cite{8}, $k=1$)

 \begin{mydef92}
\be \label{1.11}
\left|\frac
{\zeta\left(\frac 12+i\alpha_1^2\right)}
{\zeta\left(\frac 12+i\beta_1^2\right)}
\right|\sim
\frac{\sqrt{\zeta(2\sigma)}}{|\zeta[\sigma+i\alpha_0^2(\sigma)]|},\ \sigma\geq 1 .
\ee
 \end{mydef92}

\begin{mydef93}
\be \label{1.12}
\begin{split}
& \left|
\frac
{\zeta\left(\frac 12+i\alpha_1^3\right)}
{\zeta\left(\frac 12+i\beta_1^3\right)}
\right|\sim
\pi^l\sqrt{c_l}
\left|
\int_0^{\alpha_0^3(T)} \arg\zf {\rm d}t
\right|^{-l} .
\end{split}
\ee
\end{mydef93}

\begin{mydef94}
\be \label{1.13}
\begin{split}
& \left|
\frac
{\zeta\left(\frac 12+i\alpha_1^1\right)}
{\zeta\left(\frac 12+i\beta_1^1\right)}
\right|\sim
\sqrt[4]{\frac{2\pi}{H}}\frac{1}{\sqrt{|\zeta(\frac 12+i\alpha_0^1)|}}.
\end{split}
\ee
\end{mydef94}

\begin{remark}
Of course, the results (\ref{1.9}), (\ref{1.11}) -- (\ref{1.13}) are based on properties of the Jacob's ladder
$\vp_1(t)$ as follows:
\be \label{1.14}
\begin{split}
& \alpha_0^4=\vp_1^1(d)=\vp_1(d)\in (2\pi L,2\pi L+U), \\
& \alpha_1^4=\vp_1^0(d)=d\in (\overset{1}{\wideparen{2\pi L}},\overset{1}{\wideparen{2\pi L+U}}), \\
& \beta_1^4=\vp_1^0(e)=e\in (\overset{1}{\wideparen{2\pi L}},\overset{1}{\wideparen{2\pi L+U}}),
\end{split}
\ee
say (comp. (\ref{1.14}) and (\ref{1.10})).
\end{remark}

\subsection{}

Let us remind - for completeness - that Jacob's ladder
\bdis
\vp_1(t)=\frac 12\vp(t)
\edis
has been introduced in our work \cite{4} (see also \cite{5}), where the function
\bdis
\vp(t)
\edis
is arbitrary solution of the non-linear integral equation
\bdis
\int_0^{\mu[x(T)]}Z^2(t)e^{-\frac{2}{x(T)}t}{\rm d}t=\int_0^T Z^2(t){\rm d}t,
\edis
where each admissible function $\mu(y)$ generates the solution
\bdis
y=\vp(T;\mu)=\vp(T),\ \mu(y)\geq 7y\ln y.
\edis
The function $\vp_1(t)$ is called the Jacob's ladder according to Jacob's dream in Chumash, Bereishis, 28:12.

\begin{remark}
We have shown (see \cite{4}), by making use of these Jacob's ladders, that the classical Hardy-Littlewood integral (1918)
\bdis
\int_0^T\left|\zf\right|^2{\rm d}t
\edis
has - in addition to the Hardy-Littlewood expression (and other similar to that one) possessing an unbounded error at $T\to\infty$ -
the following infinite set of almost exact expressions
\bdis
\begin{split}
& \int_0^T\left|\zf\right|^2{\rm d}t = \vp_1(T)\ln\vp_1(T)+ \\
& + (c-\ln 2\pi)\vp_1(T)+c_0+\mcal{O}\left(\frac{\ln T}{T}\right),\ T\to\infty,
\end{split}
\edis
where $c$ is the Euler's constant and $c_0$ is the constant from the Titchmarsh-Kober-Atkinson formula.
\end{remark}

\begin{remark}
The Jacob's ladder $\vp_1(t)$ can be interpreted by our formula (see \cite{4})
\bdis
T-\vp_1(T)\sim (1-c)\pi(T);\ \pi(T)\sim\frac{T}{\ln T},
\edis
where $\pi(T)$ is the prime-counting function, as an asymptotically complementary function to
\bdis
(1-c)\pi(T)
\edis
in the following sense
\bdis
\vp_1(T)+(1-c)\pi(T)\sim T,\ T\to\infty .
\edis
\end{remark}

\section{Factorization, oscillating Q-system and its metamorphoses as a generic complement to the
Riemann's functional equation on the critical line}

\subsection{}

The oscillating Q-system was defined in our work \cite{7}, (2.1) as follows
\be \label{2.1}
\begin{split}
& G(x_1,\dots,x_k;y_1,\dots,y_k)\overset{def}{=}\prod_{r=1}^k
\left|
\frac{\zeta\left(\frac 12+ix_r\right)}{\zeta\left(\frac 12+iy_r\right)}
\right|= \\
& =\prod_{r=1}^k
\left|
\frac
{\sum_{n\leq\tau(x_r)}\frac{2}{\sqrt{n}}\cos\{ \vth(x_r)-x_r\ln n\}+R(x_r)}
{\sum_{n\leq\tau(y_r)}\frac{2}{\sqrt{n}}\cos\{ \vth(y_r)-y_r\ln n\}+R(y_r)}
\right|, \\
& \tau(t)=\sqrt{\frac{t}{2\pi}},\ R(t)=\mcal{O}(t^{-1/4}),\ k\leq k_0\in\mbb{N},
\end{split}
\ee
for corresponding sets (see \cite{7}, (2.2)) of the points
\bdis
(x_1,\dots,x_k),\ (y_1,\dots,y_k).
\edis

\begin{remark}
It is clear that the definition relation (\ref{2.1}) is based on simple generalization
\bdis
\left|
\frac{\zeta\left(\frac 12+ix\right)}{\zeta\left(\frac 12+iy\right)}
\right| \longrightarrow \prod_{r=1}^k
\left|
\frac{\zeta\left(\frac 12+ix_r\right)}{\zeta\left(\frac 12+iy_r\right)}
\right|
\edis
(comp. (\ref{1.7}), (\ref{1.11}) -- (\ref{1.13})).
\end{remark}

Let us remind some of the previous results playing the role of the Riemann's functional equation on the
critical line.

\subsection*{(A)}

There are the functions (see \cite{7}, (2.5))
\bdis
\begin{split}
& \alpha_r^2=\alpha_r^2(\sigma,T,\Theta,k,\epsilon),\ r=0,1,\dots,k, \\
& \beta_r^2=\beta_r^2(T,\Theta,k),\ r=1,\dots,k, \\
& \alpha_r^2,\beta_r^2\not=\gamma:\ \zeta\left(\frac 12+i\gamma\right)=0,
\end{split}
\edis
for admissible
\bdis
\sigma,T,\Theta,k,\epsilon
\edis
such that the following factorization formula
\be \label{2.2}
\prod_{r=1}^k
\left|
\frac{\zeta\left(\frac 12+i\alpha_r^2\right)}{\zeta\left(\frac 12+i\beta_r^2\right)}
\right|\sim
\frac{\sqrt{\zeta(2\sigma)}}{|\zeta[\sigma+i\alpha_0^2(\sigma)]|},\ T\to\infty,
\ee
(see \cite{7}, (2.6)) holds true, i.e. there is following set of metamorphoses of the oscillating Q-system (\ref{2.1}):
\be \label{2.3}
\begin{split}
& \prod_{r=1}^k
\left|
\frac
{\sum_{n\leq\tau(\alpha_r^2)}\frac{2}{\sqrt{n}}\cos\{ \vth(\alpha_r^2)-\alpha_r^2\ln n\}+R(\alpha_r^2)}
{\sum_{n\leq\tau(\beta_r^2)}\frac{2}{\sqrt{n}}\cos\{ \vth(\beta_r^2)-\beta_r^2\ln n\}+R(\beta_r^2)}
\right| \sim \\
& \sim \sqrt{\zeta(2\sigma)}
\left|\sum_{n=1}^\infty\frac{\mu(n)}{\sigma+i\alpha_0^2(\sigma)}\right|,\quad \sigma>1+\epsilon,\ T\to\infty
\end{split}
\ee
(see \cite{7}, (2.6)), where $\mu(n)$ is the M\" obius function.

\subsection*{(B)}

There are functions
\bdis
\begin{split}
& \alpha_r^3=\alpha_r^3(T,l,\epsilon,k),\ r=0,1,\dots,k,\ l\in\mbb{N} \\
& \beta_r^3=\beta_r^3(T,\epsilon,k),\ r=1,\dots,k, \\
& \alpha_r^3,\beta_r^3\not=\gamma:\ \zeta\left(\frac 12+i\gamma\right)=0,
\end{split}
\edis
for admissible
\bdis
T,l,\epsilon,k
\edis
such that the following factorization formula
\be \label{2.4}
\begin{split}
& \left|\int_0^{\alpha_0^3}\arg\zf {\rm d}t\right|\sim \\
& \sim \pi c^{\frac{1}{2l}}\prod_{r=1}^k
\left|
\frac{\zeta\left(\frac 12+i\alpha_r^3\right)}{\zeta\left(\frac 12+i\beta_r^3\right)}
\right|^{-\frac{1}{l}},\ T\to\infty
\end{split}
\ee
(see \cite{8}, (2.4)) holds true, i.e. there is following set of metamorphoses of the oscillating Q-system (\ref{2.1}):
\be \label{2.5}
\begin{split}
& \prod_{r=1}^k
\left|
\frac
{\sum_{n\leq\tau(\alpha_r^3)}\frac{2}{\sqrt{n}}\cos\{ \vth(\alpha_r^3)-\alpha_r^3\ln n\}+R(\alpha_r^3)}
{\sum_{n\leq\tau(\beta_r^3)}\frac{2}{\sqrt{n}}\cos\{ \vth(\beta_r^3)-\beta_r^3\ln n\}+R(\beta_r^3)}
\right| \sim \\
& \sim \pi^l \sqrt{c_l}\left|\int_0^{\alpha_0^3}\arg\zf {\rm d}t\right|,\ T\to\infty,
\end{split}
\ee
(see \cite{8}, (4.6)).

\subsection*{(B1)}

If we rewrite the formula (\ref{2.4}) as follows
\bdis
\begin{split}
 & \left|\int_{\mu_m}^{\alpha_0^3}\arg\zf{\rm d}t\right|\sim \\
 & \sim \pi c_l^{\frac{1}{2l}}\prod_{r1=1}^k
 \left|
\frac
{\sum_{n\leq\tau(\alpha_r^3)}\frac{2}{\sqrt{n}}\cos\{ \vth(\alpha_r^3)-\alpha_r^3\ln n\}+R(\alpha_r^3)}
{\sum_{n\leq\tau(\beta_r^3)}\frac{2}{\sqrt{n}}\cos\{ \vth(\beta_r^3)-\beta_r^3\ln n\}+R(\beta_r^3)}
\right|^{-\frac 1l}
\end{split}
\edis
where
\bdis
m=m(\alpha_0^3),\ \mu_m<\alpha_0^3<\mu_{m+1};\ S_1(\mu_m)=0
\edis
(see \cite{8}, (4.11)), then we obtain the set of metamorphoses (\ref{2.5}) in reverse direction: we begin with
\bdis
\left|\int_{0}^{w}\arg\zf{\rm d}t\right|
\edis
that is the Aaron staff (say),
\bdis
\longrightarrow \left|\int_{\mu_m}^{\alpha_0^3}\arg\zf{\rm d}t\right|
\edis
that is the bud of the Aaron staff (corresponding to $w=\alpha_0^3$)
\bdis
\sim \left|
\frac
{\sum_{n\leq\tau(\alpha_r^3)}\frac{2}{\sqrt{n}}\cos\{ \vth(\alpha_r^3)-\alpha_r^3\ln n\}+R(\alpha_r^3)}
{\sum_{n\leq\tau(\beta_r^3)}\frac{2}{\sqrt{n}}\cos\{ \vth(\beta_r^3)-\beta_r^3\ln n\}+R(\beta_r^3)}
\right|
\edis
already metamorphosed into almonds ripened (comp. Chumash, Bamidbar, 17:23).

\subsection*{(C)}

We have obtained the first set of metamorphoses of the primeval multiform
\bdis
G(x_1,\dots,x_k)=\prod_{r=1}^k \left|\zeta\left(\frac 12+ix_r\right)\right|
\edis
in our paper \cite{6}. The corresponding results expressed in terms of oscillating Q-system
(\ref{2.1}) are (see \cite{6}, (1.7), (2.5)): the factorization formula
\be \label{2.6}
\prod_{r=1}^k
\left|
\frac{\zeta\left(\frac 12+ix\right)}{\zeta\left(\frac 12+iy\right)}
\right|\sim \sqrt[4]{\frac{2\pi}{H}}\frac{1}{|\zeta(\frac 12+i\alpha_0^1)|}
\ee
and corresponding set of metamorphoses of the oscillating Q-system
\be \label{2.7}
\begin{split}
 & \prod_{r=1}^k
 \left|
\frac
{\sum_{n\leq\tau(\alpha_r^1)}\frac{2}{\sqrt{n}}\cos\{ \vth(\alpha_r^1)-\alpha_r^1\ln n\}+R(\alpha_r^1)}
{\sum_{n\leq\tau(\beta_r^1)}\frac{2}{\sqrt{n}}\cos\{ \vth(\beta_r^1)-\beta_r^1\ln n\}+R(\beta_r^1)}
\right| \sim \\
& \sim \frac{\sqrt[4]{\frac{2\pi}{H}}}
{\sqrt{\left|\sum_{n\leq\tau(\alpha_0^1)}\frac{2}{\sqrt{n}}\cos\{ \vth(\alpha_0^1)-\alpha_0^1\ln n\}+R(\alpha_0^1)\right|}},\
T\to\infty .
\end{split}
\ee

\subsection*{(D)} Moreover, the sequences
\bdis
\{\alpha_r^n\}_{r=0}^k,\ \{\beta_r^n\}_{r=1}^k,\ n=1,2,3
\edis
have the following universal propery:
\bdis
\begin{split}
 & \alpha_{r+1}^n-\alpha_r^n\sim (1-c)\pi(T),\ r=0,1,\dots,k-1, \\
 & \beta_{r+1}^n-\beta_r^n\sim (1-c)\pi(T),\ r=1,\dots,k-1,\ k\geq 2,
\end{split}
\edis
(comp. (\ref{1.6}), (\ref{1.7}), Remark 8 and 9).

\section{Theorem; factorization as an analogue of the Riemann's functional equation on the critical line}

\subsection{}
We use the following Euler's integral (see \cite{1}, pp. 134, 135)

\be \label{3.1}
\begin{split}
 & \int\frac{{\rm d}\vp}{a+b\cos\vp}=\frac{1}{\sqrt{a^2-b^2}}\arctan\frac{(a-b)\tan\frac{\vp}{2}}{\sqrt{a^2-b^2}},\
 \vp\in (0,\pi), \\
 & a+b>0, a^2-b^2>0 \ \Longrightarrow \ a>|b|.
\end{split}
\ee
We obtain immediately from (\ref{3.1}) the following formula
\be \label{3.2}
\begin{split}
 & \int_{2\pi L}^{2\pi L+U}\frac{{\rm d}\vp}{a+b\cos\vp}= \\
 & = \frac{1}{a+b}U
 \frac
 {\arctan\left(\sqrt{\frac{a-b}{a+b}}\tan\frac U2\right)}
 {\sqrt{\frac{a-b}{a+b}}\frac U2},\ L\in\mbb{Z},\ U\in (0,\pi).
\end{split}
\ee

\subsection{}

Now, if we use our method of transformation (see \cite{7}, (4.1) -- (4.19)) in the case of the formula (\ref{3.2}) then
we obtain the following statement.

\begin{mydef1}
 Let
\be \label{3.3}
[2\pi L,2\pi L+U]\longrightarrow [\overset{1}{\wideparen{\pi L}},\overset{1}{\wideparen{\pi L+U}}],\dots,
[\overset{k}{\wideparen{\pi L}},\overset{k}{\wideparen{\pi L+U}}]
\ee
where
\bdis
[\overset{r}{\wideparen{\pi L}},\overset{r}{\wideparen{\pi L+U}}],\ r=1,\dots,k,\ k\leq k_0\in\mbb{N}
\edis
are reversely iterated segments corresponding to the first segment in (\ref{3.3}) and $k_0$ be arbitrary and fixed number.
Then there is a sufficiently big
\bdis
T_0=T_0(a,b)>0
\edis
such that for every
\bdis
L>\frac{1}{2\pi}T_0
\edis
and for every admissible $L,U,k$ there are functions
\be \label{3.4}
\begin{split}
 & \alpha_r^4(L,U,k;a,b),\ r=0,1,\dots,k , \\
 & \beta_r^4(L,U,k), \ r=1,\dots,k, \\
 & \alpha_r^4,\beta_r^4\not=\gamma:\ \zeta\left(\frac 12+i\gamma\right)=0
\end{split}
\ee
such that
\be \label{3.5}
\begin{split}
 & \prod_{r=1}^k
 \left|
\frac{\zeta\left(\frac 12+i\alpha_r^4\right)}{\zeta\left(\frac 12+i\beta_r^4\right)}
\right|^2\sim \\
& \sim
\frac
 {\arctan\left(\sqrt{\frac{a-b}{a+b}}\tan\frac U2\right)}
 {\sqrt{\frac{a-b}{a+b}}\frac U2}
 \frac{a+b\cos\alpha_0^4}{a+b},\ L\to \infty,
\end{split}
\ee
of course
\bdis
G(\alpha^4,\beta^4)\sim \sqrt{(\dots)} \ \Longleftrightarrow \
\{ G(\alpha^4,\beta^4)\}^2\sim (\dots).
\edis
Moreover, the sequences
\bdis
\{\alpha_r^4\}_{r=0}^k,\ \{\beta_r^4\}_{r=1}^k,\ n=1,2,3
\edis
have the following properties
\be \label{3.6}
\begin{split}
 & 2\pi L<\alpha_0^4<\alpha_1^4<\dots<\alpha_k^4, \\
 & 2\pi L<\beta_1^4<\beta_2^4<\dots<\beta_k^4, \\
 & \alpha_0^4\in (2\pi L,2\pi L+U),\\
 & \alpha_r^4,\beta_r^4\in (\overset{r}{\wideparen{2\pi L}},\overset{r}{\wideparen{2\pi L+U}}),
  r=1,2,\dots,k,
\end{split}
\ee
\be \label{3.7}
\begin{split}
 & \alpha_{r+1}^4-\alpha_r^4\sim (1-c)\boldmath{\pi}(2\pi L),\ r=0,1,\dots,k-1, \\
 & \beta_{r+1}^4-\beta_r^4\sim (1-c)\boldmath{\pi}(2\pi L),\ r=1,\dots,k-1,\ k\geq 2,
\end{split}
\ee
where
\bdis
\boldmath{\pi}(T)\sim \frac{T}{\ln T},\ T\to\infty
\edis
is the prime-counting function and $c$ is the Euler's constant.
\end{mydef1}

\subsection{}

Now, let us notice the following.

\begin{remark}
The asymptotic behavior of the following sets
\be \label{3.8}
\{\alpha_r^4\}_{r=0}^k
\ee
is as follows: if $L\to\infty$ then the points of every set (\ref{3.8}) recede unboundedly each from other and all
together recede to infinity. Hence, at $L\to\infty$ each of the sets (\ref{3.8}) behaves as one-dimensional
Friedmann-Hubble universe.
\end{remark}

\begin{remark}
Next, we express the result (\ref{3.5}) in connection with (\ref{1.1}), (\ref{1.3}), Remark 2 and Remark 9 in the form
\be \label{3.9}
\prod_{r=0}^k\left|\zeta\left(\frac 12+i\alpha_r^4\right)\right|\sim
\sqrt{\chi_4(U,\alpha_0^4)}\prod_{r=1}^k\left|\zeta\left(\frac 12+i\beta_r^4\right)\right|,
\ee
where $\chi_4$ stands for the right-hand side of (\ref{3.5}). It is quite clear that this formula if an analogue of the
Riemann's function equation on the critical line.
\end{remark}

\begin{remark}
Of course, the first result (\ref{1.9}) is a particular case of (\ref{3.5}).
\end{remark}

\section{On a set of metamorphoses that correspond to the formula (\ref{3.5})}

\subsection{}

Let us remind the spectral form of the Riemann-Siegel formula
\be \label{4.1}
\begin{split}
 & Z(t)=\sum_{n\leq \tau(x_r)}\frac{2}{\sqrt{n}}\cos\left\{ t\omega_n(x_r)+\psi(x_r)\right\}+
 \mcal{O}(x_r^{-1/4}), \\
 & \tau(x_r)=\sqrt{\frac{x_r}{2\pi}}, \\
 & t\in [x_r,x_r+V],\ V\in (0,x_r^{1/4}) ,
\end{split}
\ee
and similarly for $x_r\longrightarrow y_r$, where
\bdis
T_0<2\pi L<x_r,y_r
\edis
(see \cite{6}, (6.1), comp. \cite{8}, (4.4) and Remark 6 ibid).

\begin{remark}
We call the expressions
\be \label{4.2}
\frac{2}{\sqrt{n}}\cos\left\{ t\omega_n(x_r)+\psi(x_r)\right\}, \dots
\ee
as the Riemann's oscillators with:
\begin{itemize}
 \item[(a)] the amplitude
 \bdis
 \frac{2}{\sqrt{n}},
 \edis
 \item[(b)] the incoherent local phase-constant
 \bdis
 \psi(x_r)=-\frac{x_r}{2}-\frac{\pi}{8},
 \edis
 \item[(c)] the non-synchronized local time
 \bdis
 t=t(x_r)\in [x_r,x_r+V],
 \edis
 \item[(d)] the local spectrum of the cyclic frequencies
 \bdis
 \{\omega_n(x_r)\}_{n\leq \tau(x_r)},\ \omega_n(x_r)=\ln\frac{\tau(x_r)}{n},\ \dots
 \edis
\end{itemize}
\end{remark}

\begin{remark}
Now, we see that the Q-system, where of course
\be \label{4.3}
\left|\zf\right|=|Z(t)|,
\ee
is really complicated oscillating system (comp. (\ref{2.1}), (\ref{4.1}), (\ref{4.3})).
\end{remark}

\subsection{}

Consequently, we have (see (\ref{2.1}), (\ref{3.5}), (\ref{4.1}), Remark 14 and (\ref{4.3})) the following.

\begin{mydef41}
The following set of metamorphoses
\be \label{4.4}
\begin{split}
 & \prod_{r=1}^k
 \left|
\frac
{\sum_{n\leq \tau(\alpha_r^4)}\frac{2}{\sqrt{n}}\cos\{ \alpha_r^4\omega_n(\alpha_r^4)+\psi(\alpha_r^4)\}+R(\alpha_r^4)}
{\sum_{n\leq \tau(\beta_r^4)}\frac{2}{\sqrt{n}}\cos\{ \beta_r^4\omega_n(\beta_r^4)+\psi(\beta_r^4)\}+R(\beta_r^4)}
 \right|\sim \\
 & \sim
\frac
{\arctan\left(\sqrt{\frac{a-b}{a+b}}\tan\frac U2\right)}
{\sqrt{\frac{a-b}{a+b}}\frac U2}
\frac{a+b\cos(\alpha_0^4)}{a+b},\ L\to\infty.
\end{split}
\ee
correponds to the factorization formula (\ref{3.5}).
\end{mydef41}

\subsection{}

Now, let us notice the following.

\begin{remark}
By Theorem, there are control functions (\ref{3.4}) (Golem's shem) of the set of metamorphoses (\ref{4.4}) of the
oscillating Q-system (\ref{2.1}), (see also (\ref{4.1}), (\ref{4.3})).
\end{remark}

\begin{remark}
The mechanism of metamorphosis is as follows. Let (comp. (\ref{3.4}) and \cite{7}, (2.2))
\be \label{4.5}
\begin{split}
 & M_k^3=\{ \alpha_1^4,\dots,\alpha_k^4\}, \\
 & M_k^4=\{ \beta_1^4,\dots,\beta_k^4\},
\end{split}
\ee
where. of course, (comp. \cite{7}, (2.12))
\be \label{4.6}
\begin{split}
 & M_k^3\subset M_k^1\subset (T_0,+\infty)^k, \\
 & M_k^4\subset M_k^2\subset (T_0,+\infty)^k .
\end{split}
\ee
Now, if we obtain, after random sampling such points (comp. conditions \cite{7}, (2.2)) that
\be \label{4.7}
\begin{split}
 & (x_1,\dots,x_k)=(\alpha_1^4,\dots,\alpha_k^4)\subset M_k^3, \\
 & (y_1,\dots,y_k)=(\beta_1^4,\dots,\beta_k^4)\subset M_k^4,
\end{split}
\ee
(see (\ref{4.5}), (\ref{4.6})), then -- at the points (\ref{4.7}) -- the Q-system (\ref{2.1}) changes its old form
(=chrysalis) into its new form (=butterfly) and the last is controlled by the function $\alpha_0^4$.
\end{remark}

\begin{remark}
Now, it should be clear that the set of metamorphoses of oscillating Q-system also belongs to the family of
analogues of the Riemann's functional equation on the critical line.
\end{remark}

\section{On decomposition of the result of metamorphoses (\ref{4.4}) into three parts: signal, noise and error term}

In this section we use the terminology from the theory of signal processing.

\subsection{}

Let us remind (see \cite{7}, (4.11)) that
\bdis
\begin{split}
 & \tilde{Z}^2(t)=\frac{\left|\zf\right|^2}{\omega(t)}, \\
 & \omega(t)=\left\{ 1+\mcal{O}\left(\frac{\ln\ln t}{\ln t}\right)\right\}\ln t.
\end{split}
\edis
Since in our case
\bdis
t\longrightarrow 2\pi L,
\edis
then (comp. \cite{7}, (4.11), (4.12))
\bdis
\tilde{Z}^2(\alpha_r^4)=
\frac{\left|\zeta\left(\frac 12+i\alpha_r^4\right)\right|^2}
{\left\{ 1+\mcal{O}\left(\frac{\ln\ln L}{\ln L}\right)\right\}\ln t},\ \dots
\edis

\begin{remark}
Consequently, the primary form of the asymptotic formula for metamorphoses (\ref{4.4}) is as follows
\be \label{5.1}
\begin{split}
 & \prod_{r=1}^k
 \left|
\frac
{\sum_{n\leq \tau(\alpha_r^4)}\frac{2}{\sqrt{n}}\cos\{ \alpha_r^4\omega_n(\alpha_r^4)+\psi(\alpha_r^4)\}+R(\alpha_r^4)}
{\sum_{n\leq \tau(\beta_r^4)}\frac{2}{\sqrt{n}}\cos\{ \beta_r^4\omega_n(\beta_r^4)+\psi(\beta_r^4)\}+R(\beta_r^4)}
 \right|\sim \\
 & \sim \left\{ 1+\mcal{O}\left(\frac{\ln\ln L}{\ln L}\right)\right\}
 \frac
{\arctan\left(\sqrt{\frac{a-b}{a+b}}\tan\frac U2\right)}
{\sqrt{\frac{a-b}{a+b}}\frac U2}
\frac{a+b\cos(\alpha_0^4)}{a+b} .
\end{split}
\ee
\end{remark}

Since the last two factors on the right-hand side of (\ref{5.1}) are bounded functions for all
\bdis
U\in (0,\pi),\ L>\frac{1}{2\pi} T_0
\edis
(for all fixed admissible $k,a,b$, see (\ref{3.1}), (\ref{3.4})), then we obtain from (\ref{5.1}) the following.

\begin{mydef42}
\be \label{5.2}
\begin{split}
 & \prod_{r=1}^k
 \left|
\frac
{\sum_{n\leq \tau(\alpha_r^4)}\frac{2}{\sqrt{n}}\cos\{ \alpha_r^4\omega_n(\alpha_r^4)+\psi(\alpha_r^4)\}+R(\alpha_r^4)}
{\sum_{n\leq \tau(\beta_r^4)}\frac{2}{\sqrt{n}}\cos\{ \beta_r^4\omega_n(\beta_r^4)+\psi(\beta_r^4)\}+R(\beta_r^4)}
 \right|= \\
 & = \frac{a}{a+b}
 \frac
{\arctan\left(\sqrt{\frac{a-b}{a+b}}\tan\frac U2\right)}
{\sqrt{\frac{a-b}{a+b}}\frac U2}+ \\
& + \frac{b}{a+b}\frac
{\arctan\left(\sqrt{\frac{a-b}{a+b}}\tan\frac U2\right)}
{\sqrt{\frac{a-b}{a+b}}\frac U2}\cos(\alpha_0^4)+ \\
& + \mcal{O}\left(\frac{\ln\ln L}{\ln L}\right),\ L\to\infty .
\end{split}
\ee
\end{mydef42}

\subsection{}

Let us remind (see (\ref{3.4})), that
\be \label{5.3}
\alpha_0^4=\alpha_0^4(L,U,k;a,b)=\alpha_0^4(L,U)
\ee
for admissible and fixed $k,a,b$.

\begin{itemize}
 \item[(a)] We see that the first function on the right-hand side of (\ref{5.2}) is the $L$-th member
\be \label{5.4}
\begin{split}
 & f(2\pi L+U)=g_L(U)=\frac{a}{a+b}
 \frac
{\arctan\left(\sqrt{\frac{a-b}{a+b}}\tan\frac U2\right)}
{\sqrt{\frac{a-b}{a+b}}\frac U2},\ U\in (0,\pi) \\
& g_L(U)=g_{L'}(U),\ \forall\ L,L'>\frac{1}{2\pi}T, \\
& 2\pi L+U\in [2\pi L,2\pi L+\pi)
\end{split}
\ee
of the stationary sequence
\be \label{5.5}
\{ g_L(U)\}_{L>T_0/2\pi},\ U\in (0,\pi).
\ee
By (\ref{5.4}), (\ref{5.5}) the corresponding signal is defined.

Consequently, by the first function on the right-hand side of (\ref{5.2}) deterministic signal is
expressed (see (\ref{5.5}).
\item[(b)] The main factor in the second member is the following function
\bdis
\cos(\alpha_0^4),\ \alpha_0^4\in \alpha_0^4(L,U),
\edis
where (comp. (\ref{1.14}), (\ref{5.3}))
\bdis
\alpha_0^4=\vp_1(d)\in (2\pi L,2\pi L+U),\ d=d(L,U),
\edis
and $\vp_1(d)$ is the value of the Jacob's ladder. That is, the distribution of the values
\bdis
\alpha_0^4\in (2\pi L,2\pi L+U),\ U\in (0,\pi)
\edis
we may suppose as very complicated.

Consequently, the second function we shall characterize as noise -- a non-useful part of the signal. The noise
may be controlled by variation of the parameter $b$,
\bdis
a>|b|
\edis
(see (\ref{3.1})), i.e. by abatement of $|b|$.
\item[(c)] The third function we shall call (fine) error term, since
\bdis
\mcal{O}\left(\frac{\ln\ln L}{\ln L}\right)\xrightarrow{L\to\infty}0.
\edis
\end{itemize}

\begin{remark}
Hence, the final state of metamorphoses in (\ref{5.2}) is split into three parts: signal, noise and (fine) error term.
\end{remark}

\section{The set of distinct metamorphoses in (\ref{4.4})}

Of course, there is a point $U_0\in (0,\pi)$ such that (see (\ref{5.4})
\bdis
g_L'(U)|_{U=U_0}\not=0,\ \forall\ L>\frac{T_0}{2\pi},
\edis
i.e.
\be \label{6.1}
\begin{split}
 & g'_L(U)\not=0,\ U\in O_\delta(U_0)=(U_0-\delta,U_0+\delta), \\
 & \forall \ U',U''\in O_\delta(U_0),\ U'\not= U'' \ \Rightarrow \ g_L(U')\not=g_L(U'')
\end{split}
\ee
for suitable $\delta>0$.

Next, we shall suppose that there are such
\bdis
U_1,U_2\in O_\delta(U_0),\ U_1\not=U_2,
\edis
that
\be \label{6.2}
\begin{split}
& \alpha_r^4(U_1,L)=\alpha_r^4(U_2,L),\ r=0,1,\dots,k, \\
& \beta_r^4(U_1,L)=\beta_r^4(U_2,L),\ r=1,\dots,k
\end{split}
\ee
for all
\bdis
L>\tilde{L}>\frac{T_0}{2\pi}.
\edis
In this case we obtain, by comparison of the formulas (\ref{5.1}), for $U_1, U_2$ that
\bdis
\begin{split}
& \frac
{\arctan\left(\sqrt{\frac{a-b}{a+b}}\tan\frac{U_1}{2}\right)}
{\sqrt{\frac{a-b}{a+b}}\frac{U_1}{2}}= \\
& = \left\{ 1+\mcal{O}\left(\frac{\ln\ln L}{\ln L}\right)\right\}
\frac
{\arctan\left(\sqrt{\frac{a-b}{a+b}}\tan\frac{U_2}{2}\right)}
{\sqrt{\frac{a-b}{a+b}}\frac{U_2}{2}},
\end{split}
\edis
i.e. in the limit case we obtain the equality
\bdis
\frac
{\arctan\left(\sqrt{\frac{a-b}{a+b}}\tan\frac{U_1}{2}\right)}
{\sqrt{\frac{a-b}{a+b}}\frac{U_1}{2}}=\frac
{\arctan\left(\sqrt{\frac{a-b}{a+b}}\tan\frac{U_2}{2}\right)}
{\sqrt{\frac{a-b}{a+b}}\frac{U_2}{2}}
\edis
that contradicts with (\ref{6.1}).

Hence, we have that for every
\bdis
U',U''\in O_\delta(U_0),\ U'\not= U''
\edis
there is an infinite subsequence
\bdis
\{ \bar{L}\}\subset \{ L\},\ \bar{L}>\tilde{L}
\edis
such that (comp. (\ref{6.2}))
\be \label{6.3}
\begin{split}
& \left( \alpha_0^4(U',\bar{L}), \alpha_1^4(U',\bar{L}),\dots,\alpha_k^4(U',\bar{L}), \right. \\
& \left. \beta_1^4(U',\bar{L}),\dots,\beta_k^4(U',\bar{L})\right) \not= \\
& \not=\left(
\alpha_0^4(U'',\bar{L}), \alpha_1^4(U'',\bar{L}),\dots,\alpha_k^4(U'',\bar{L}), \right. \\
& \left. \beta_1^4(U'',\bar{L}),\dots,\beta_k^4(U'',\bar{L})
\right),\ \forall \bar{L}\in \{ \bar{L}\}.
\end{split}
\ee
Consequently, by (\ref{6.3}) we have the following.

\begin{mydef43}
There is an infinite set of distinct metamorphoses (\ref{4.4}).
\end{mydef43}

\thanks{I would like to thank Michal Demetrian for his help with electronic version of this paper.}


\begin{thebibliography}{29}
\bibitem{1}
Leonardo Eulero, \emph{Institutionum calculi integralis}, volumen primum, Petropoli 1768; Moscow, 1956 (in Russian).
%
\bibitem{2}
G.H. Hardy, \emph{Divergent series}, Oxford, Oxford Univ. Press, 1949.
%
\bibitem{3}
E. Landau, \emph{Handbuch der Lehre von der Verteilung der Primzahlen}, Leipzig und Berlin, Teubner, 1909.
%
\bibitem{4}
J. Moser,
`Jacob's ladders and almost exact asymptotic representation of the Hardy-Littlewood integral`,
Math. Notes 88, (2010), 414-422, arXiv: 0901.3937.
%
\bibitem{5}
J. Moser,
`Jacob's ladders, structure of the Hardy-Littlewood integral and some new class of nonlinear integral equations`,
Proc. Steklov Inst. 276 (2011), 208-221, arXiv: 1103.0359.
%
\bibitem{6}
J. Moser, `Jacob's ladders, $\zeta$-factorization and infinite set of metamorphosis of a multiform`, arXiv: 1501.07705.
%
\bibitem{7}
J. Moser,
`Jacob's ladders, Riemann's oscillators, quotient of the oscillating multiforms and set of metamorphoses of this system`,
arXiv: 1506.00442.
%
\bibitem{8}
J. Moser,
`Jacob's ladders, new properties of the function $\arg\zf$ and corresponding metamorphoses`,
arXiv: 1506.07967.
%
\bibitem{9}
E.C. Titchmarsh, `\emph{The theory of the Riemann zeta-function}`, Clarendon Press, Oxford, 1951.



\end{thebibliography}
\end{document}